# TOWARDS OPTIMAL ENERGY-WATER SUPPLY SYSTEM OPERATION FOR AGRICULTURAL AND METROPOLITAN ECOSYSTEMS

Marcello Di Martino[1, 2], Patrick Linke[1, 3] and Efstratios N. Pistikopoulos[1, 2, *]
1. Artie McFerrin Department of Chemical Engineering, Texas A&M University, College Station, TX, USA
2. Texas A&M Energy Institute, Texas A&M University, College Station, TX, USA
3. Department of Chemical Engineering, Texas A&M University at Qatar, Education City, Doha, Qatar

*Abstract*

The energy-water demands of metropolitan regions and agricultural ecosystems are ever-increasing. To tackle this challenge efficiently and sustainably, the interdependence of these interconnected resources has to be considered. In this work, we present a holistic decision-making framework which takes into account simultaneously a water and energy supply system with the capability of satisfying metropolitan and agricultural resource demands. The framework features: (i) a generic large-scale planning and scheduling optimization model to minimize the annualized cost of the design and operation of the energy-water supply system, (ii) a mixed-integer linear optimization formulation, which relies on the development of surrogate models based on feedforward artificial neural networks and first-order Taylor expansions, and (iii) constraints for land and water utilization enabling multi-objective optimization. The framework provides the operational profiles of all energy-water system elements over a given time horizon, which uncover potential synergies between the essential food, energy, and water resource supply systems.

*Keywords*

Food-Energy-Water Nexus, Surrogate Modeling, Mixed-Integer Optimization, Multi-objective Optimization.

**Introduction**

Metropolitan and agricultural ecosystems are facing increased water and energy demands due to rising populations. In addition, societies are trying to improve their standard of living through economic growth, while natural resources are degrading. Overall, this puts immense stresses on the resource supply systems (Daher and Mohtar, 2015). In order to capture the dynamics, operability and complexity of these systems, surrogate modeling and optimization are of the utmost importance to evaluate trade-offs between the cost and the resource utilization of supply system alternatives. To consider sustainability, the interdependence of food, energy, and water resources must be modeled, naturally resulting in the incorporation of a food-energy-water nexus (FEWN) methodology approach (Garcia and You, 2016), where systematic and integrated analysis of process structure alternatives is enabled through superstructure synthesis (Di Martino et al., 2020, 2021). However, optimizing the design of interconnected process systems and evaluating interactions and trade-offs between sub-systems remains challenging (Ibrić et al., 2021).

To uncover trade-offs between technology alternatives and synergies between interconnected resource systems, a time-dependent integrated decision-making framework is necessary (Namany et al., 2019). Here, we present a holistic planning and scheduling optimization framework incorporating the energy and water supply system for metropolitan and agricultural resource demand satisfaction. The energy supply system based on the selection of available renewable energy and energy storage technologies is optimized, with the objective to minimize total annual cost over a predefined time horizon. Operational parameters of a three-stage reverse osmosis (RO) desalination plant are determined by supplying the energy for the RO system exclusively through the defined renewable energy technologies. Vice versa, the water demands of the energy supply system must be satisfied by the RO plant. Thus, the framework consists of a detailed and interconnected energy-water supply system. Therefore, periods of high renewable resource availability can be used to produce a surplus of desalinated water for later usage in biomass farming or water demand satisfaction, instead of increasing energy storage, since this option might result in a cost advantage. Accordingly, the main contributions of this work are the linkage of a renewable energy and a RO water supply system design and operation optimization framework, as well as the consideration of desalinated water storage as a form of energy storage.

The remainder of this paper is structured as follows: The problem statement is introduced together with an overview of the considered system. Afterwards, the utilized surrogate models are presented, together with the overall optimization model. Then, the capabilities of the derived optimization framework are illustrated in a case study.

**Energy-Water Nexus – Problem Statement**

This work applies a planning and scheduling optimization formulation to an energy-water supply system to simultaneously satisfy specified energy and water demand profiles over a predefined time horizon. The optimization framework can not only be used to satisfy the water and energy needs of a metropolitan area, but also to supply the necessary energy-water resources for farming purposes. Therefore, the proposed energy-water supply decision-making framework encompasses an energy and water sub-system, which is illustrated in Figure 1. The energy sub-system is based on the work of Cook et al. (2022) and consists of the selected energy sources, energy conversion technologies and energy storage technologies. To retain the linearity of the model, linear energy conversion technology surrogate models have to be derived, supplying correlations for cost and land use of each technology dependent on the power output. Furthermore, the energy profiles of all considered renewable energies must be defined. For this work, we consider wind farms, solar farms based on single axis tracking and fixed angle, as well as biomass in the form of maize as viable renewable energy technology options. To consider the water supply sub-system, a RO desalination plant model is utilized according to the work of Di Martino et al. (2022), since RO is the industry leading desalination technology. The RO model specifies all operational pressures together with the permeate output, to derive the energy consumption of the plant. It is important to note, that the necessary water for potential biomass farming is supplied by the RO system, whereas the energy demands of the RO system must be satisfied by the energy system. Overall, obtained solutions specify the cost minimal design and operation of the energy-water supply system, meaning over the predefined time horizon hourly energy production and energy storage levels are specified together with the RO plant's water recovery, permeate output and energy consumption. Additionally, ε-constraints for the total land use and total water use are introduced to enable multi-objective optimization. No losses of the water and energy storage systems are considered.

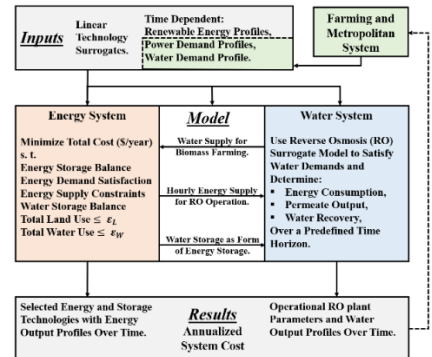

*Figure 1. Schematic overview of the proposed energy-water supply decision-making framework for agricultural and metropolitan applications.*

This planning and scheduling model results in a large scale mixed-integer programming problem. To maintain solvability, a linear model structure is enforced by introducing linear surrogates for all observed nonlinearities (Pochet and Wolsey, 2006, Mahdi and Dawood, 2022). Overall, this approach yields a mixed-integer linear programming (MILP) problem. The model is solved in GAMS using the Gurobi solver.

**Model Development**

Surrogate models based on either first-order Taylor expansions or piece-wise linear approximations are utilized to approximate nonlinearities. All surrogates have been derived in MATLAB.

*Energy Supply Sub-System*

According to the work of Cook et al. (2022) linear technology surrogates to derive the necessary land and annualized cost of each energy conversion technology dependent on the annual power output can be derived. To consider biomass farming in the form of maize all relevant information could be found in the literature (Braun et al., 2010, Lunik and Langemeier, 2015). To obtain linear correlations for solar and wind farms, optimization models were utilized to minimize the cost of the system while satisfying a predefined power demand and deriving the corresponding land use. These mixed-integer linear optimization models generally have the same form as summarized in Eq. (1) to (9). Here, $t \in T$ denotes the time horizon and time scale under investigation, whereas $n \in \bar{N}$ represents the location and numbering of selected energy technology units, i.e., wind turbines or solar panels.

$$\min C^{Total} = C^{Land} + C^{Tech} \tag{1}$$
$$s.t. \ C^{Tech} = k_{Tech} \cdot N \tag{2}$$
$$C^{Land} = k_{Land} \cdot A \cdot N \tag{3}$$
$$P^{Sum} \geq P^{Target} \tag{4}$$
$$P^{Sum} = \sum_{n,t} p_{n,t} \cdot y_{n,t}^{op} \tag{5}$$
$$y_{n,t}^{op} \leq y_n^{buy} \quad \forall t \in T, n \in \bar{N} \tag{6}$$
$$N = \sum_n y_n^{buy} \tag{7}$$
$$p_{n,t} = f_1(p_{in,t}, n) \quad \forall t \in T, n \in \bar{N} \tag{8}$$
$$A = f_2(A^{Tech}, A^{Spacing}) \tag{9}$$

The total cost $C^{Total}$ [\$/year] is composed of the land cost $C^{Land}$ and the technology cost $C^{Tech}$, as stated in Eq. (1). $C^{Tech}$ can be derived based on the number of units $N$ and the annualized technology cost factor $k_{Tech}$ (Eq. (2)). $C^{Land}$ can similarly be derived based on the annualized land cost $k_{Land}$, necessary area per unit $A$ [ha] and $N$ (Eq. (3)). It is important to note that both $k_{Tech}$ and $k_{Land}$ are dependent on the location of the technology under investigation. Equation (4) ensures that the desired power target ($P^{Target}$ [kWh/year]) is fulfilled by the generated power of the system ($P^{Sum}$). $P^{Sum}$ is derived based on the power output of a single unit over time and location $p_{n,t}$ [kW] when it is in operation (Eq. (5)). Operation and existence of units is monitored with the binary variables $y_{n,t}^{op}$ and $y_n^{buy}$, which are correlated as shown in Eq. (6). Accordingly, $N$ is derived based on $y_n^{buy}$ in Eq. (7). Equation (8) states that the power output of a single unit is calculated technology dependent (considered with function $f_1$) based on the location of the unit $n$ and renewable power input $p_{in,t}$, i.e., solar DNI or wind speed. Lastly, Eq. (9) denotes the technology dependent calculation (considered with function $f_2$) of $A$ based on the area of a unit itself $A^{Tech}$ and the required spacing between units $A^{Spacing}$.

Then, $P^{Target}$ can be varied within a relevant range to obtain land use and system cost dependent on the desired power output. The obtained points are approximated with a linear regression based on minimizing the squared error. This approach results in overall eight model surrogates, with mean-square errors ($R^2$) between 0.97 and 1. Regarding energy storage technologies, information concerning the expected lifespan, efficiency, as well as capital and operational cost per capacity need to be supplied.

*Water Supply Sub-System*

As presented in Di Martino et al. (2022), a feedforward artificial neural network (ANN) with rectified linear units (ReLU) as activation functions is used to approximate the permeate concentration and the pressure difference across the energy recovery device, since it is possible to reformulate these ANNs as MILPs. Further, linear correlations based on the feed pressures of each stage and each parallel flow are used to calculate the retentate pressure of each stage and parallel flow, as well as the water recovery of each stage. Remaining nonlinearities consist of the calculation of the overall system water recovery ($WR^{sys}$) based on the water recovery of each stage ($WR_1, WR_2, WR_3$), the calculation of the permeate volume flow ($Q_p$) based on the feed volume flow ($Q_f$) and $WR^{sys}$, as well as the calculation of the energy consumption ($EC$) of the system. To approximate $WR^{sys}$ and $Q_p$ a first-order Taylor expansion around the nominal operating point ($\overline{WR^{sys}} = 0.6039, \overline{WR_1} = 0.3113, \overline{WR_2} = 0.2935, \overline{WR_3} = 0.1860, \overline{Q_p} = 975 \frac{m^3}{h}$), obtained from one-and-a-half years of hourly RO operational data, is performed as shown from Eq. (10) to (15).

$$WR^{sys} = WR_1 + (1 - WR_1) \cdot WR_2 + (1 - WR_1) \cdot (1 - WR_2) \cdot WR_3 \tag{10}$$

$$\nabla WR^{sys} = \begin{pmatrix} 1 - WR_2 - WR_3 + WR_2 \cdot WR_3 \\ 1 - WR_1 - WR_3 + WR_1 \cdot WR_3 \\ 1 - WR_1 - WR_2 + WR_1 \cdot WR_2 \end{pmatrix} \tag{11}$$

$$WR^{sys} \approx \overline{WR^{sys}} + \nabla WR^{sys}(\overline{WR_1}, \overline{WR_2}, \overline{WR_3})^T \cdot \begin{pmatrix} WR_1 - \overline{WR_1} \\ WR_2 - \overline{WR_2} \\ WR_3 - \overline{WR_3} \end{pmatrix} \tag{12}$$

$$Q_p = WR^{sys} \cdot Q_f \tag{13}$$

$$\nabla Q_p = \begin{pmatrix} Q_f \\ WR^{sys} \end{pmatrix} \tag{14}$$

$$Q_p \approx \overline{Q_p} + \frac{\overline{Q_p}}{\overline{WR^{sys}}} \cdot (WR^{sys} - \overline{WR^{sys}}) + \overline{WR^{sys}} \cdot (Q_f - \frac{\overline{Q_p}}{\overline{WR^{sys}}}) \tag{15}$$

For the relevant ranges of $0.4 \leq WR^{sys} \leq 0.85$ and 227 m³/h $\leq Q_p$ the given correlations are sufficiently accurate approximations with $R^2$ values of 0.98 for both cases. Lastly, the energy consumption ($EC$ in [kW]) of the RO system is approximated with a piece-wise linear approximation. Originally, the non-linear specific energy consumption of the RO plant could be approximated with a linear reformulation (Di Martino et al., 2022). However, to obtain from the former term the overall energy consumption, the specific energy consumption has to be multiplied with $WR^{sys}$ and $Q_f$. This results in a trilinear term, which could be approximated to an accuracy of $R^2 = 0.99$ with a shallow ANN with one hidden layer and two nodes with ReLUs as activation functions. The ANN is trained with the aforementioned one-and-a-half years of operational RO data and subsequently reformulated as a MILP. Thus, ultimately resulting in a piecewise linear approximation of the energy consumption. The inputs of the ANN have been identified as $WR_1, WR_2, WR_3$ and $Q_f$. More information regarding the RO desalination operational plant data can be found in Di Martino et al. (2022). Details on how to reformulate an ANN with ReLUs as activation function can be found in Grimstad and

Andersson (2019). The obtained weights and biases of the ANN approximating $EC$ are summarized in Table 1. Table 2 gives an overview and comparison of the employed surrogate models.

*Table 1. Used weights and biases of the shallow ANN (1 hidden layer, two nodes) to approximate $EC$ with $WR_1, WR_2, WR_3$ and $Q_f$.*

| Layer | Weight | Bias |
|---|---|---|
| Hidden Layer | 1.1556 -0.5436 0.4116 0.4630 0.2843 0.1137 0.1071 1.3843 | 1.0465 -0.3829 |
| Output Layer | 0.4638 0.5371 | -0.9950 |

*Table 2. Summary and comparisons of the presented surrogate models.*

| Method | Input | Output | $R^2$ |
|---|---|---|---|
| Taylor Expansion | $WR_1, WR_2, WR_3$ | $WR^{sys}$ | 0.98 |
| Taylor Expansion | $WR^{sys}, Q_f$ | $Q_p$ | 0.98 |
| ANN with ReLUs | $WR_1, WR_2, WR_3, Q_f$ | $EC$ | 0.99 |

**Energy-Water Optimization Model**

The decision-making optimization framework consists of the presented first-order Taylor expansions for $WR^{sys}(t)$ and $Q_p(t)$, the MILP reformulation for $EC(t)$ based on $WR_1(t), WR_2(t), WR_3(t)$ and $Q_f(t)$, the "combined MILP model" as described in Cook et al. (2022) and the RO desalination plant model as summarized in Di Martino et al. (2022). The referenced RO system is modified to a time dependent model by introducing a time index for all variables enabling the investigation of time dependent interconnected energy-water supply system operation. To enforce resource utilization, $WR^{sys}(t)$ is restricted according to Eq. (16).

$$WR^{sys}(t) \geq WR^{Lim} \quad \forall t \in T \tag{16}$$

Furthermore, the water balance equations connecting the water and energy supply system for all $t \in T$ are stated from Eq. (17) to (21). Here, the linear $Q_p$ surrogate enables the direct usage of the RO permeate output in other constraints instead of the indirect approach of restricting $Q_f$ and $WR_1, WR_2, WR_3$ as done in Di Martino et al. (2022).

$$Q_p(t) = Q_{En}^D(t) + Q_W^D(t) - Q_{stor}(t) + Q_{rel}(t) \tag{17}$$
$$V(t) = V(t-1) + Q_{stor}(t) \cdot \Delta t - Q_{rel}(t) \cdot \Delta t \tag{18}$$
$$Q_{rel}(t_1) \cdot \Delta t \leq V(t_1) \tag{19}$$
$$V(t_1) \leq V(t_{|T|}) \tag{20}$$
$$V_T \geq V(t) \tag{21}$$

In Eq. (17) $Q_{En}^D(t)$ and $Q_W^D(t)$ denote the water demand of the energy system and all other applications, i.e., metropolitan and agricultural systems, respectively. Further, $Q_{stor}(t)$ is the amount of water stored and $Q_{rel}(t)$ the amount of water released from the storage system. The water storage balances presented form Eq. (18) to Eq. (21) introduces the overall tank volume $V_T$ and the water storage level at each time point $V(t)$, as well as beginning of time and end of time constraints to ensure operability and applicability of obtained solutions.

The energy balance equations for all $t \in T$ from Cook et al. (2022) are modified according to Eq. (22) to (25).

$$P^{Sum} = \sum_t P(t) - EC(t) + \sum_{t,k}[P_{rel,k}(t) - P_{stor,k}(t)] \tag{22}$$
$$P_T^D \leq P^{Sum} \tag{23}$$
$$P(t) - EC(t) + \sum_k[P_{rel,k} - P_{stor,k}(t)] \geq P^D(t) \tag{24}$$
$$P(t) = \sum_{Tech} n_{Tech} \cdot p_{Tech}(t) \tag{25}$$

Here, the novelty is to introduce the energy consumption of the RO system in the overall energy balance (Eq. (22)) and energy calculation at each time point (Eq. (24)). Accordingly, energy demand targets to be satisfied in each time point ($P^D(t)$) and over the complete time horizon ($P_T^D$) can be introduced. The energy derived from all utilized renewable energy technologies (summarized in set $Tech =$ {single axis tracking, fixed angle solar panels, wind turbines, maize farming}) is calculated as shown in Eq. (25). In this case, $n_{Tech}$ denotes the number of units of each renewable energy technology and $p_{Tech}(t)$ the energy output of a single technology unit based on the observed climate data.

Thus, the option of building up a desalinated water storage instead of energy storage during periods of high renewable energy availability is enabled through the incorporation of $EC(t)$ in Eq. (22) and (24). The amount of water stored or released from storage based on either a surplus or deficit of produced permeate is derived with Eq. (17). The water storage level at each time point is calculated according to Eq. (18), where for $t =1$ an initial water storage level is permitted ($V(0)$). The tank sizing equations are summarized from Eq. (19) to (21).

The objective of the framework is to minimize the annual cost of the energy-water supply system. Accordingly, the objective function presented in Cook et al. (2022) is extended by the investment and operational cost of the RO system, together with the cost of water storage. The investment and operational RO cost based on the RO plant's capacity $Q_p^{Cap}$ are shown in Eq. (26) and (27), with an assumed operating life of 20 years (Wittholz et al., 2008). $Q_p^{Cap}$ is in turn derived according to the observed $Q_p(t)$ values, as illustrated in Eq. (28). The cost of water storage dependent on the tank volume is displayed in Eq. (29) with an assumed lifetime of 30 years (Loh et al., 2002).

$$C_{Inv}^{RO} = 1.25 \cdot 10^{-5} \cdot Q_p^{Cap} + 10.429 \cdot 10^6 \quad (26)$$
$$C_{Op}^{RO} = 0.45 \cdot \sum_t Q_p(t) \quad (27)$$
$$Q_p^{Cap} \geq Q_p(t) \quad \forall \, t \in T \quad (28)$$
$$C^T = 0.75 \cdot V_T + 5000 \quad (29)$$

The size of the resulting optimization model is dependent on the selected timescale and time horizon (Brunaud et al., 2019). The size of the model for hourly evaluation is illustrated for varying selected time horizons in Table 3. The goal is to solve this model for a one-year time horizon at an hourly time scale ($|T| = 8760$). However, in this case the size of the model is computationally intractable without advanced algorithmic solution approaches like Benders or Dantzig-Wolfe decomposition strategies. Furthermore, the framework is intended to be used for various scenario analyses and multi-objective optimization. Therefore, fast solution generation is essential (within 1000s to 10000s). To overcome this challenge, without using advanced algorithmic strategies, we propose employing a steady-state assumption for the water supply sub-system. Ultimately, this means that one optimal operational point of the RO desalination plant is specified over one year, while the energy sub-system is solved at an hourly time scale. The obtained solution acts as an upper bound of the true optimum and can be used as an initial cost estimator. In this case, solutions can be obtained within 120s.

*Table 3. Size of the energy-water decision making model dependent on the time horizon under investigation.*

| Timeframe | #Constraints | #Continuous Variables | #Binaries |
|---|---|---|---|
| 1 day | 2596 | 2306 | 265 |
| 2 weeks | 35980 | 32278 | 3697 |
| 1 month | 77068 | 69142 | 7921 |
| 6 months | 468688 | 420502 | 48181 |
| 1 year | 937348 | 840982 | 96361 |

**Case Study**

To illustrate the developed decision-making framework and underline its capabilities, it is applied to an energy-water supply case study in Texas. Here, renewable energies in the form of wind speed and solar DNI are plentiful available. Also, the used RO data is based on a desalination plant in this location. The objective of this study is to evaluate changing resource restrictions in terms of water and land utilization on the solution. Specifically, how the optimal energy mix and operating points change by adjusting $\varepsilon_L$ and $\varepsilon_W$ (see Figure 1). The supply system is utilized to satisfy the water and energy demand of a metropolitan region combined with greenhouse farming. It is assumed that a water demand of 570 m³/h must be satisfied at all times by the RO system. Furthermore, the water recovery must be at least 60% for efficient resource utilization. Additionally, the energy demand of 100 kW together with the energy and water demand of 1 greenhouse must be satisfied at all time points. Yara International, Qatar Fertiliser Company and Hassad Food farmed tomatoes for a one-year time period in Qatar in a pilot water saving greenhouse. The results have been summarized in a trial report which specifies the greenhouse's demand profiles, which are applied to this case study due to comparable climate conditions between the two locations. Initially, the water and land use of the supply system are not restricted. Then, both resources are simultaneously restricted by successively reducing $\varepsilon_L$ and $\varepsilon_W$. By adjusting these values in different ranges other potential energy and water applications can be considered. The results of this study are summarized in Figure 2, which illustrates how the annual energy demand is supplied for various water and land restrictions, together with the system cost. The RO system operates at $WR^{sys} = 60\%$ with $Q_P$ between 568 m³/h and 600 m³/h, while the $SEC$ is approx. constant at 0.3 kWh/m³.

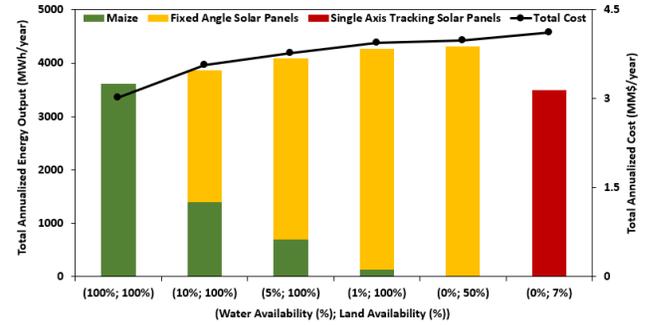

*Figure 2. Annualized energy output by technology (bars) and annualized cost (points) of energy-water supply systems for varying water and land availability.*

To underline the framework's versatility, an energy base load demand of 1 GW together with the energy and water demand of 10 greenhouses needs to be satisfied at all times in a subsequent study. In this case, for all investigated land and water restrictions, wind turbines supply at least 99% of the energy, resulting in a system with annualized costs between 1.45 and $2.32 \cdot 10^9$ \$/year and total energy outputs between 0.899 and $1.07 \cdot 10^3$ GWh/year. The RO system only constitutes at most 0.35% of the annual cost, while the wind farm and the energy storage system denote at least 99% of the total cost. For the unrestricted case, the wind farm and the storage system contributed 73.5% and 26% of the total cost, whereas in the most restricted case these fractions change to 38.9% and 61%, respectively.

Interestingly, in both investigated case studies, for all presented solutions, water storage as a form of energy storage is not selected. This is also the case when the investment cost of water storage is neglected in the objective function. Intuitively, water storage should be the cheaper option to energy storage, however due to the steady-state approximation of the water-system and the possibility of too high specific energy consumptions of the RO system this might not be the case. In this model, energy storage is not connected to energy loss, whereas the production of water results in energy losses to the system due to the energy demands of the RO process. To further analyze this observation, the complete time dependent system is solved for a two-week time horizon at an hourly time scale for an energy demand of 100 kW in addition to one greenhouse, together with a 1% water use restriction and no restriction on land use. In this case we obtain a solution with a water storage tank level as displayed in Figure 3. Accordingly, water storage should be considered as an alternative to energy storage in coupled dynamic energy-water supply systems. The system cost in this case is $87,000 with a total energy output of 590,000 kWh/year supplied by one wind turbine.

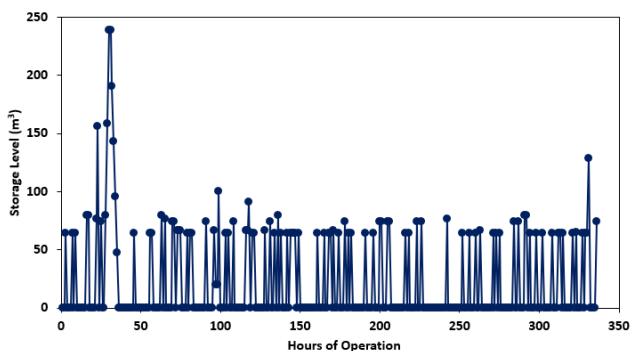

*Figure 3. Water tank storage level for two weeks of hourly operation (336 hours).*

## Conclusion

We presented an integrated planning and scheduling energy-water supply system model capable of satisfying varying energy and water demand profiles, dependent on the desired application, for minimizing the system's annual cost by deriving the design and all operational parameters of the system. The capabilities of the framework have been illustrated with a case study focusing on the design and operation of an energy-water supply system in Texas based on solely renewable energy supply. The proposed system can be used to satisfy the water and energy demands of metropolitan and agricultural systems, thus incorporating all pillars of the FEWN. Moreover, multi-objective analysis has been performed by evaluating water and land resource restrictions. This in turn underlines the sensitivity of solutions to available resource restrictions. It is important to note that the presented model is generic in nature and can be applied to any region and scale of interest by changing the supplied model inputs.

## Acknowledgments

This publication was made possible by the National Priorities Research Program (NPRP) grant No. NPRP11S-0107–180216 from the Qatar National Research Fund (a member of Qatar Foundation). The findings herein reflect the work, and are solely the responsibility, of the authors. The authors also gratefully acknowledge support from Texas A&M University, the Texas A&M Energy Institute and Texas A&M University at Qatar.